\documentclass{amsart}

\usepackage{epsfig}
\usepackage{amscd}
\usepackage{amsmath, amssymb, amsthm}
\usepackage{graphics}
\usepackage{color}
\usepackage{dsfont}
\newtheorem*{main-theorem}{Main Theorem}
\newtheorem{proposition}{Proposition}[section]
\newtheorem{theorem}{Theorem}
\newtheorem{lemma}[proposition]{Lemma}
\newtheorem{corollary}[proposition]{Corollary}

\theoremstyle{definition}

\newtheorem{remark}[proposition]{Remark}

\numberwithin{equation}{section}

\def\11{\mathbf{1}}

\def\reals{{\mathbb R}}
\def\cx{{\mathbb C}}

\def\Ci{{\mathcal C}^\infty}

\def\Im{\,\mathrm{Im}\,}

\def\O{{\mathcal O}}
\def\SS{{\mathbb S}}

\def\phi{\varphi}
\def\half{{\frac{1}{2}}}
\def\dist{\text{dist}\,}

\def\be{\begin{eqnarray*}}
\def\ee{\end{eqnarray*}}
\def\ben{\begin{eqnarray}}
\def\een{\end{eqnarray}}

\def\lll{\left\langle}
\def\rrr{\right\rangle}

\def\L2R{L_{\text{Rest}}^2}

\def\L2c{L^2_{\text{comp}}}

\def\CC{\mathcal{C}}
\def\vol{\text{vol}}

\begin{document}
\title[Resolvents and SLS]{Cutoff Resolvent Estimates and the
  Semilinear Schr\"odinger Equation}
\author{Hans Christianson}
\address{Department of Mathematics, University of California, Berkeley, CA 94720 USA}
\email{hans@math.berkeley.edu}

\begin{abstract}
This paper shows how abstract resolvent estimates imply local
smoothing for solutions to the Schr\"odinger equation.  If the
resolvent estimate has a loss when compared to the optimal,
non-trapping estimate, there is a corresponding loss in regularity in
the local smoothing estimate.  As an application, we apply well-known
techniques to obtain well-posedness results for the semi-linear
Schr\"odinger equation.
\end{abstract}
\maketitle


\section{Introduction}

In this short note we show how cutoff semiclassical resolvent
estimates for the Laplacian on a non-compact manifold, with spectral
parameter on the real axis, lead to well-posedness results for the
semilinear 
Schr\"odinger equation.  Motivated by the requirements of \cite{Ch3}
and \cite{BGT}, and the microlocal inverse estimates of \cite{Ch,
  Ch2}, we first prove a general theorem for a large class of
resolvents.  Following the recent work of Nonnenmacher-Zworski
\cite{NZ}, we apply the general theorem in the case there is a
hyperbolic fractal trapped set.  

Let $(M,g)$ be a Riemannian manifold of dimension $n$ without
boundary, with
(non-negative) 
Laplace-Beltrami operator $- \Delta$ acting on functions.  The Laplace-Beltrami
operator is an unbounded, essentially self-adjoint operator on
$L^2(M)$ with domain $H^2(M)$.  We assume $(M,g)$ is asymptotically
Euclidean in the sense of \cite[(3.7)-(3.9)]{NZ} and that the classical resolvent $(-\Delta - (\lambda^2 +i \epsilon))^{-1}$ obeys
a limiting absorption principle as $ \epsilon
\to 0+$, $\lambda \neq 0$.

Our first result is
that if we have cutoff semiclassical resolvent estimates with a
sufficiently small loss, then we have weighted smoothing for the
Schr\"odinger propagator with a loss.  Let $\rho_s$ be a smooth,
non-vanishing weight function satisfying
\ben
\label{rho-def}
\rho_s(x) \equiv \lll d_g(x, x_0 )\rrr^{-s},
\een
for some fixed $x_0$ and $x$ outside a compact set.

\begin{theorem}
\label{Sch-thm}
Suppose for each compactly supported function $\chi \in \Ci_c(M)$ with
sufficiently small support, there is $h_0>0$ such that the semi-classical Laplace-Beltrami operator satisfies
\ben
\label{abs-res-est}
\|\chi ( - h^2\Delta - E )^{-1} \chi u \|_{L^2(M)} \leq \frac{g(h)}{h} \|  u
\|_{L^2(M)}, \,\,\, E >0
\een
uniformly in $0<h \leq h_0$, where $g(h) \geq c_0 >0$, $g(h) =
o(h^{-1})$.  Then for each $T>0$ and $s > 1/2$, there is a
constant $C=C_{T,s}>0$
such that 
\ben
\label{Sch-thm-est-2}
\int_0^T \left\| \rho_s e^{it \Delta  } u_0 \right\|_{H^{1/2 -
    \eta}(M)}^2 dt  \leq C \| u_0 \|_{L^2(M)}^2,
\een
where $\eta \geq 0$ satisfies
\ben
\label{eta-cond}
g(h) h^{2 \eta } = \O(1),
\een
and $\rho_s$ is given by \eqref{rho-def}.
\end{theorem}



The assumption that $(M,g)$ is asymptotically Euclidean is that there
exists $R_0>0$ sufficiently large that, on
each infinite branch of $M \setminus B(0,R_0)$, the semiclassical
Laplacian $-h^2 \Delta$ takes the form
\be
-h^2 \Delta|_{M \setminus B(0,R_0)} = \sum_{|\alpha| \leq 2}
a_\alpha(x,h) (hD_x)^\alpha,
\ee
with $a_\alpha(x,h)$ independent of $h$ for $|\alpha|=2$, 
\be
&& \sum_{|\alpha| = 2}
 a_\alpha(x,h) (hD_x)^\alpha \geq C^{-1} |\xi|^2, \,\,\, 0 < C <
\infty, \text{ and}\\
&& \sum_{|\alpha| \leq 2}
a_\alpha(x,h) (hD_x)^\alpha \to |\xi|^2, \,\,\, \text{as } |x| \to
\infty \text{ uniformly in } h.
\ee
In order to quote the results of \cite{NZ} we also need the following
analyticity assumption: $\exists \theta_0 \in [0, \pi)$ such that the
$a_\alpha(x,h)$ are extend holomorphically to
\be
\{ r \omega : \omega \in \cx^n, \,\, \dist ( \omega, \SS^n) <
\epsilon, \,\, r \in \cx, \,\, |r| \geq R_0, \,\, \arg r \in
[-\epsilon, \theta_0+ \epsilon) \}.
\ee
As in \cite{NZ}, the analyticity assumption immediately implies 
\be
\partial_x^\beta \left( \sum_{|\alpha| \leq 2} a_\alpha(x, h)
  \xi^\alpha - |\xi|^2 \right)  = o ( |x|^{-|\beta|} )\lll \xi \rrr^2,
\,\, |x| \to \infty.
\ee


Recall the free Laplacian $(-\Delta_0 - \lambda^2 )^{-1}$ on $\reals^n$
has a holomorphic continuation from $\Im \lambda >0$ to $\lambda \in
\cx$ for $n \geq 3$ odd, and to the logarithmic covering space for $n$ even.  This motivates the limiting absorption assumption, that
\be
\lim_{\epsilon \to 0+, \,\,\, \lambda \neq 0} \rho_s (- \Delta -
(\lambda^2 +i\epsilon))^{-1} \rho_s
\ee
exists as a bounded operator 
\be
L^2(M, d \vol_g)
\to L^2(M, d \vol_g),
\ee
provided $s >1/2$.  As in the free case, we allow a possible logarithmic singularity at $\lambda = 0$.

The problem of ``local smoothing'' estimates for the Schr\"odinger
equation has a long history.  The sharpest results to date are those
of Doi \cite{Doi} and Burq \cite{Bur2}.  Doi proved if $M$ is
asymptotically Euclidean, then one has the estimate
\ben
\label{doi-est}
\int_0^T \left\| \chi e^{it\Delta} u_0 \right\|_{H^{1/2}(M)}^2 dt  \leq C \| u_0 \|_{L^2(M)}^2
\een
for $\chi \in \Ci_c(M)$ if and only if there are no trapped sets.
Burq's paper showed if there is trapping due to the presence
of several convex obstacles in $\reals^n$ satisfying certain assumptions, then one has the estimate \eqref{doi-est}
with the $H^{1/2}$ norm replaced by $H^{1/2 - \eta}$ for $\eta
>0$.  In \cite{Ch3}, the author considered an arbitary, single trapped
hyperbolic orbit.  One of the goals of this paper is to use estimates
obtained by Nonnenmacher-Zworski \cite{NZ} for fractal hyperbolic
trapped sets to obtain similar results to \cite{Ch3} for the
semilinear Schr\"odinger equation.  To that end we have the following
corollary to Theorem \ref{Sch-thm}.    

\begin{corollary}
\label{Sch-cor}
Assume $(M,g)$ admits a hyperbolic fractal trapped set, $K_E$, in the energy level
$E>0$ and that the topological pressure $P_E(1/2)<0$.  Then $-h^2 \Delta -E$
satisfies \eqref{abs-res-est} for some $E>0$ with $g(h) = C \log(1/h)$, and for
every $\eta>0$, $T>0$, and $s>1/2$, there exists a constant $C =
C_{P_E,\eta,T, s}>0$ such that
\be
\int_0^T \left\| \rho_s e^{it \Delta  } u_0 \right\|_{H^{1/2 -
    \eta}(M)}^2 dt  \leq C \| u_0 \|_{L^2(M)}^2.
\ee
\end{corollary}
We remark that the assumption $P_E(1/2)<0$ implies the trapped set
$K_E$ is filamentary or ``thin'' (see \cite{NZ} for
definitions).

We consider the following semilinear Schr\"odinger equation problem:
\ben
\label{nls}
\left\{ \begin{array}{cc}
i \partial_t u + \Delta u = F(u) \,\, \text{on } I \times M; \\
u(0,x) = u_0(x),
\end{array} \right.
\een
where $I \subset \reals$ is an interval containing $0$.  Here the
nonlinearity $F$ satisfies
\be
F(u) = G'(|u|^2)u,
\ee
and $G: \reals \to \reals$ is at least $C^3$ and satisfies
\be
| G^{(k)}(r) | \leq C_k \langle r \rangle^{\beta - k },
\ee
for some $\beta \geq \half$.

In \S \ref{strichartz} we prove a family of Strichartz-type estimates which will result in the following 
well-posedness theorem.  

\begin{theorem}
\label{nls-lwp}
Suppose $(M,g)$ satisfies the assumptions of the introduction, and
set
\ben
\label{delta-def}
\delta = \frac{4 \eta}{2 \eta +1} \geq 0.
\een
Then for each 
\ben
\label{s-cond}
s > \frac{n}{2} - \frac{2}{\max \{2\beta -2, 2 \} } + \delta
\een
and each $u_0 \in H^s(M)$ there exists $p > \max \{ 2\beta -2, 2 \}$ and $0< T \leq
1$ such that \eqref{nls} has a unique solution
\ben
\label{u-soln}
u \in C ([-T, T]; H^s(M)) \cap L^p([-T, T]; L^\infty(M)).
\een

Moreover, the map $u_0(x) \mapsto u(t,x) \in C([-T,T];H^s(M))$ is Lipschitz continuous on
bounded sets of $H^s(M)$, and if $\|u_0\|_{H^s}$ is bounded, $T$ is bounded
from below.

If, in addition, $(M,g)$ satisfies the assumptions of Corollary \ref{Sch-cor}, $n \leq 3$, $\beta <3$, and $G(r) \to +
\infty$ as $r \to + \infty$, then $u$ in \eqref{u-soln} extends to a solution
\be
u \in C ((-\infty, \infty); H^1(M)) \cap L^p((-\infty, \infty);
L^\infty(M)). 
\ee
\end{theorem}
\begin{remark}
In particular, the cubic defocusing non-linear Schr\"odinger equation
is globally $H^1$-well-posed in three dimensions with a fractal trapped
hyperbolic set which is sufficiently filamentary.  Of course other
nonlinearities can be considered, but for simplicity we consider only
these in this work.
\end{remark}

{\bf Acknowledgments.} 
This
research was partially conducted during the period the author was
employed by the Clay Mathematics Institute as a Liftoff Fellow.



\section{Proof of Theorem \ref{Sch-thm}}

\begin{figure}
\centerline{\input{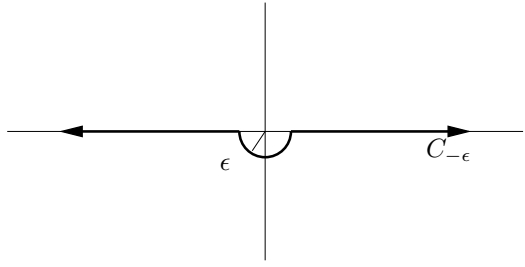}}
\caption{\label{spec101-fig} The curve $\CC_{-\epsilon}$ in the complex plane.}
\end{figure}

Since we are assuming $(-\Delta - z )^{-1}$ obeys a limiting
absorption principle, we have
\be
\| \rho_s   (-\Delta - (\tau - i
\epsilon ) )^{-1} \rho_s \|_{L^2 \to L^2} \leq C_\epsilon
\ee
for $0 < \epsilon_0 \leq |\tau| \leq C$.  
For $| \sigma| \geq C$ for some $C>0$, $\sigma \in \cx$ in a
neighbourhood of the real axis, write
\be
-\Delta - \sigma & = & - \Delta - \frac{z}{h^2} \\
& = & h^{-2} ( - h^2 \Delta - z),
\ee
for 
\be
z \in [E - \alpha, E + \alpha] + i [-c_0 h, c_0 h].
\ee
Now 
\be
 ( - h^2 \Delta - z)
\ee
is a Fredholm operator for $z$ in the specified range, and hence the
``gluing'' techniques from \cite{Vod} and \cite[\S 2]{Ch3} can be used to conclude for $s>1/2$,
\be
\rho_s ( -h^2 \Delta -z)^{-1} \rho_s
\ee
has a holomorphic extension to a slightly smaller neighbourhood in
$z$, and in particular,
\be
\| \rho_s (-h^2 \Delta -E)^{-1} \rho_s \|_{L^2 \to L^2} \leq C
\frac{g(h)}{h}.
\ee
Rescaling, we have
\ben
\label{res-tau-est}
\left\| \rho_s (-\Delta - \tau )^{-1} \rho_s
\right\|_{L^2 \to L^2 } \leq
C \frac{ g( \lll \tau \rrr^{1/2}) }{\lll \tau \rrr^{1/2}}, \,\,\, \tau
\in \CC_{\pm \epsilon},
\een
where (see Figure \ref{spec101-fig})
\be
\CC_{\pm \epsilon} = \{ \tau \in \reals : | \tau | \geq \epsilon \} \cup
\{ \tau \in \cx : | \tau | = \epsilon, \, \pm \Im \tau \geq 0 \}.
\ee

As in \cite{Ch3} and \cite{Bur2}, the following lemma follows from
integration by parts and interpolation, together with the condition 
on $\eta$, \eqref{eta-cond}.
\begin{lemma}
\label{H1-lemma}
With the notation and assumptions above, we have
\be
\| \rho_s (-\Delta  - \tau )^{-1} \rho_s \|_{L^2 \to
  H^1} \leq C g( \lll \tau \rrr^{1/2}), \,\,\, \tau \in \CC_{\pm \epsilon},
\ee
and for every $r \in [-1,1]$, 
\be
\| \rho_s (-\Delta   - \tau )^{-1} \rho_s \|_{H^r
  \to H^{1 + r - \eta/2}} \leq C, \,\,\, \tau \in  \CC_{\pm \epsilon}.
\ee
\end{lemma}

Theorem \ref{Sch-thm} now follows from the standard ``$T T^*$''
argument, letting $\epsilon \to 0$ in \eqref{res-tau-est}
(see \cite{BGT}, the references cited therein, and \cite{Ch3}).

\qed

The following Corollary uses interpolation with an $H^2$ estimate to replace the $H^{1/2- \eta}$ norm
on the left hand side of \eqref{Sch-thm-est-2} with $H^{1/2}$, and will be
of use in \S \ref{strichartz}.  See \cite{Ch3} for the details of the proof.

\begin{corollary}
\label{inter-cor-10}
Suppose $(M,g)$ satisfies the assumptions of Theorem \ref{Sch-thm}.
For each $T >0$ and $s > 1/2$, there is a constant $C>0$  such
that
\ben
\label{eqn-1002}
\int_0^T \left\| \rho_s e^{it \Delta } u_0 \right\|_{H^{1/2}(M)}^2 dt  \leq C \| u_0 \|_{H^\delta(M)}^2,
\een
where $\delta \geq 0$ is given by \eqref{delta-def}.

In particular, if $(M,g)$ satisfies the assumptions of Corollary 
\ref{Sch-cor}, then for any $\delta>0$, there is $C = C_\delta >0$
such that \eqref{eqn-1002} holds.
\end{corollary}



\section{Strichartz-type Inequalities}
\label{strichartz}
In this section we give several families of Strichartz-type
inequalities and prove Theorem \ref{nls-lwp}.  The statements and
proofs are mostly adaptations of
similar inequalities 
in \cite{BGT}, so we leave out the proofs of these in the interest of space.

If we view $M \setminus U$, where $U$ is a neighbourhood of $K_E$, as a manifold with non-trapping geometry,
we may apply the results of \cite{HTW} or \cite{BT} to a solution of the
Schr\"odinger equation away from the trapping region, resulting in
perfect Strichartz estimates.  For this section we need
\eqref{Sch-thm-est-2} only with a compact cutoff $\chi$ instead of
with the more general weight $\rho_s$.
\begin{proposition}
\label{str-prop-1}
For every $0 < T \leq 1$ and each $\chi \in \Ci_c(M)$ satisfying $\chi
\equiv 1$ near $U$ , there is a constant $C>0$ such that
\ben
\label{str-1}
\| (1-\chi) u \|_{L^p([0,T]) W^{s,q}(M)} \leq C \| u_0 \|_{H^s(M)},
\een
where $u = e^{it \Delta}u_0$, $s \in [0,1]$, and $(p,q)$, $p >2$
satisfy
\be
\frac{2}{p} + \frac{n}{q} = \frac{n}{2}.
\ee
\end{proposition}

\begin{remark}
In the sequel, wherever unambiguous, we will write
\be
L_T^pW^{s,q}:= L^p([0,T]) W^{s,q}(M)
\ee
and
\be
H^s: = H^s(M).
\ee
\end{remark}

\begin{proposition}
\label{str-prop-5}
Suppose $(M,g)$ satisfies the assumptions of the Introduction, $u =
e^{it \Delta} u_0$, and
\be
v = \int_0^t e^{i(t - \tau) \Delta} f( \tau ) d \tau.
\ee
Then for each $0 <T \leq 1$ and $\delta \geq 0$ satisfying \eqref{delta-def}, we have the estimates
\ben
\label{str-3}
\| u \|_{L^p_T W^{s- \delta, q}} \leq C \| u_0 \|_{H^s}
\een
and
\ben
\label{str-4}
\| v \|_{L_T^p W^{s - \delta, q}} \leq C \| f \|_{L_T^1 H^s},
\een
where $s\in [0,1]$ and $(p,q)$, $p>2$ satisfy
the Euclidean scaling 
\ben
\label{pq-eqn-3}
\frac{2}{p} + \frac{n}{q} = \frac{n}{2}.
\een
\end{proposition}
The proof uses a local WKB expansion localized also in time to the
scale of inverse frequency, followed by summing over frequency bands
(see \cite{Ch3} and \cite{BGT1}).  The only difference here is the
explicit dependence of $\delta$ on $\eta$, which is related to the growth of
the function $g(h)$.

\begin{proof}[Proof of Theorem \ref{nls-lwp}]
The proof of Theorem \ref{nls-lwp} is a slight modification of the
proof of Proposition 3.1 in \cite{BGT1}, but we include it here in the
interest of completeness.    
Fix $s$ satisfying \ref{s-cond} and choose $p > \max \{ 2\beta-2, 2 \}$
satisfying
\be
s > \frac{n}{2} - \frac{2}{p} + \delta \geq \frac{n}{2} - \frac{1}{\max
  \{ 2\beta-2, 2 \}}
\ee
where $\delta \geq0$ satisfies \eqref{delta-def}.  Set $\sigma = s - \delta$ and 
\be
Y_T = C ( [-T,T]; H^s(M)) \cap L^p([-T,T]; W^{\sigma,q}(M))
\ee
for 
\be
\frac{2}{p} + \frac{n}{q} = \frac{n}{2},
\ee
equipped with the norm
\be
\| u \|_{Y_T} = \max_{|t| \leq T} \| u(t)\|_{H^s(M)} + \|u\|_{L^p_T
  W^{\sigma,q}}.
\ee
Let $\Phi$ be the nonlinear functional
\be
\Phi ( u ) = e^{it \Delta} u_0 - i \int_0^t e^{i(t-\tau) \Delta}
F(u(\tau)) d \tau.
\ee
If we can show that $\Phi: Y_T \to Y_T$ and is a contraction on a ball
in $Y_T$ centered at $0$ for sufficiently small $T>0$, this will
prove the first assertion of the Proposition, along with the Sobolev
embedding
\ben
\label{sob-emb-1a}
W^{\sigma,q}(M) \subset L^\infty(M),
\een
since $\sigma >n/q$.  From Proposition \ref{str-prop-5}, we bound the
$W^{\sigma}$ part of the $Y_T$ norm by the $H^s$ norm, giving
\be
\| \Phi(u) \|_{Y_T} & \leq & C \left( \| u_0 \|_{H^s} + \int_{-T}^T
  \| F(u (\tau)) \|_{H^s} d \tau \right)\\
& \leq & C \left( \| u_0 \|_{H^s} + \int_{-T}^T \|(1 +
  |u(\tau)|)\|_{L^\infty}^{2\beta-2} ) \|u( \tau ) \|_{H^s} d \tau
\right),
\ee
where the last inequality follows by our assumptions on the structure
of $F$.  
Applying H\"older's inequality in time with $\tilde p = p/ (2\beta-2)$
and $\tilde q$ satisfying
\be
\frac{1}{\tilde q} + \frac{1}{ \tilde p} = 1
\ee
gives
\be
\| \Phi(u) \|_{Y_T}  \leq C \left( \| u_0 \|_{H^s} + T^\gamma  \|
  u\|_{L^\infty_T H^s}\|( 1 + |u| )\|_{L^p_T L^\infty}^{2\beta-2})
\right)
\ee
where $\gamma = 1/ \tilde q >0$.  Thus
\be
\| \Phi(u) \|_{Y_T} \leq  C \left( \| u_0 \|_{H^s} +
  T^\gamma(\|u\|_{Y_T} + \|u\|_{Y_T}^{ 2\beta} ) \right).
\ee
Similarly, we have for $u, v \in Y_T$,
\ben
\label{phi-ineq-2} \lefteqn{ \| \Phi(u) - \Phi(v) \|_{Y_T}  \leq } \\
 & \leq &  C T^\gamma \| u - v
\|_{L^\infty_T H^s} \|( 1 + | u|) \|_{L^p_T L^\infty}^{2\beta-2} +
\|(1 + |v|) 
\|_{L^p_T L^\infty}^{2\beta-2}) \\
& \leq & C T^\gamma \| u - v \|_{Y_T} \|( 1 + |u|) \|_{Y_T}^{2\beta-2} +
\|(1+ |v|) \|_{Y_T}^{2\beta-2}), \nonumber
\een
which is a contraction for sufficiently small $T$.  This concludes the
proof of the first assertion in the Proposition.  

To get the second assertion, we observe from \ref{phi-ineq-2} and the
definition of $Y_T$, if $u$ and $v$ are two solutions to \eqref{nls}
with initial data $u_0$ and $u_1$ respectively, so 
\be
\widetilde{\Phi}(v) =  e^{it \Delta} u_1 - i \int_0^t e^{i(t-\tau) \Delta}
F(v(\tau)) d \tau,
\ee
 we have
\be
\lefteqn{\max_{|t| \leq T } \| u(t) - v(t) \|_{H^s} } \\
& = & \max_{|t| \leq T } \| \Phi(u)(t) - \widetilde{\Phi} (v)(t) \|_{H^s} \\
& \leq C &
\Bigg( \| u_0 - u_1 \|_{H^s} \\
&& \quad + T^\gamma \max_{|t| \leq T } \| u(t) -
  v(t) \|_{H^s} \|( 1 + | u|) \|_{L^p_T L^\infty}^{2\beta-2} + \| (1
  + |v|)
  \|_{L^p_T L^\infty}^{2\beta-2} ) \Bigg),
\ee
which, for $T>0$ sufficiently small gives the Lipschitz continuity.

If $(M,g)$ satisfies the assumptions of Corollary \ref{Sch-cor}, $n \leq 3$, $\beta < 3$, and $G(r) \to + \infty$ as $r \to +\infty$, we can take
$s$ and $p$ satisfying $p > \max \{ 2 \beta -2, 2 \}$ and 
\be
s > \frac{n}{2} - \frac{2}{p} + \delta \geq \frac{n}{2} - \frac{2}{\max
  \{ 2\beta-2, 2 \}}
\ee
for any $\delta>0$.  Then $\sigma = s- \delta > q/n$ and the preceding
argument holds.  Finally, the proof of the global well-posedness now follows from the standard
global well-posedness arguments from, for example, \cite[Chapter 6]{Caz}.
\end{proof}


\end{document}